\documentclass{article}
\usepackage{hyperref}
\usepackage{subcaption}
\usepackage[utf8]{inputenc}

\usepackage{soul}

\usepackage{wrapfig}

\usepackage{comment}
\usepackage{xcolor}
\usepackage[colorinlistoftodos]{todonotes}

\newcommand{\rachel}[1]{\todo[size=\small,inline,color=violet!30]{#1 \\ \hfill --- Rachel}}
\definecolor{munsell}{rgb}{0.0, 0.5, 0.69}

\definecolor{green1}{rgb}{0.0, 0.9, 0.0}
\newcommand{\jose}[1]{\todo[size=\small,inline,color=green1]{\sf #1 \\ \hfill --- Jos\'e}}

\title{Quantifying and Documenting Inequity in PhD-granting Mathematical Sciences Departments in the United States}

\author{Ron Buckmire\footnote{Corresponding author: ron@oxy.edu} \and Carrie Diaz Eaton
\and Joseph E. Hibdon, Jr. \and Jakini Kauba \and Drew Lewis  \and Omayra Ortega \and Jos\'e L. Pab\'on \and Rachel Roca \and Andr\'es R. Vindas-Mel\'endez}

\date{\today}

\begin{document}
\maketitle

\begin{abstract}
We provide an example of the application of quantitative techniques, tools, and topics from mathematics and data science to analyze the mathematics community itself in order to quantify and document inequity in our discipline. 
This work is a contribution to the new and growing interdisciplinary field recently termed ``mathematics of Mathematics,'' or ``MetaMath.''
Using data about PhD-granting institutions in the United States and publicly available funding data from the National Science Foundation, we highlight inequalities in departments at U.S. institutions of higher education that produce PhDs in the mathematical sciences.
Specifically, we determine that a small fraction of mathematical sciences departments receive a large majority of federal funding awarded to support mathematics in the United States.
Additionally, we identify the extent to which women faculty members are underrepresented in  mathematical sciences PhD-granting institutions in the United States.
We also show that this underrepresentation of women faculty is even more pronounced in departments that received more federal grant funding. 
\end{abstract}

\newcommand{\prestige}{Wapman-prestige}
\section{Introduction}
\label{sec:intro}

The kinds of problems mathematics and data science can be used to solve are extremely varied, running the gamut from theoretical with no foreseen applications to those that are immediately applicable to important real-world phenomena like climate change, epidemiology, and social networks. 
There is a rapidly growing body of work using tools from the mathematical sciences to analyze the  mathematical sciences itself that has recently been described as the ``mathematics of Mathematics'' or ``MetaMath'' \cite{MathofMathSurvey}. 
This term was chosen as an allusion to the broader field of ``science of Science'' \cite{fortunato2018science} or ``SciSci" that uses scientific and mathematical tools to analyze science as a whole as well as its individual disciplines.

In this paper, we present a contribution to the mathematics of Mathematics to introduce readers to the kinds of questions that can be asked and the types of tools and techniques that can be used in the emerging MetaMath area. 
We are inspired by the work of Wapman, Zhang, Larremore, and Clauset \cite{wapman2022quantifying} published in \textit{Nature} in 2022 that quantified hierarchy in faculty hiring and retention in a wide array of academic disciplines in the United States by analyzing a very large dataset of nearly 300,000 faculty members in over 10,000 departments at almost 400 PhD-granting institutions from 2011 to 2020.
The research we present here builds on and leverages the ideas and analyses presented in Wapman et al. \cite{wapman2022quantifying}, but focuses on specific disciplines in the mathematical sciences that are available in their dataset \cite{wapman-data}, namely, mathematics, statistics, and operations research.
This results in a dataset which includes nearly 10,000 faculty members in well over 200 departments that granted PhDs in the mathematical sciences from 2011 to 2020.

Using mathematical tools from network science \cite{kawakatsu2021emergence}, Wapman et al. produced a ``prestige ranking'' for every department at every PhD-granting institution in their dataset and quantified hierarchy in faculty hiring and retention in a large number of academic disciplines.
In particular, Wapman et al. constructed a directed graph of PhD-granting departments with an edge from department A to department B corresponding to a faculty member who earned a Ph.D. at department A being hired into department B. 
Unlike other rankings like those published by \textit{U.S. News \& World Report} that purport to assign prestige based on an arbitrary selection of attributes, Wapman et al.'s prestige rankings are based upon the characteristics of the directed graph of departments. It is assumed  that ``prestigious" departments will prefer to hire faculty from other ``prestigious" departments, and then the overall hierarchy of departments is inferred from the topology of the network using the SpringRank algorithm \cite{de2018physical}.   
We will refer to the  definition of prestige from Wapman et al. \cite{wapman2022quantifying} as ``Wapman-prestige." 

We utilize Wapman's quantification of prestige to produce a definition of ``elite," which we will  use to help us quantify and document inequity in PhD-granting departments in the United States.
We  define ``elite" institutions as those that are within the top quartile as defined by \prestige.
In this paper we will use Hasty et al.'s anthropological definition of inequity. They say inequity is ``the unequal distribution of resources due to an unjust power imbalance. It is a type of inequality caused by this unequal distribution, often as a result of injustices against historically excluded groups of people" \cite{hasty2022anthro}.
In the mathematical sciences, the underrepresentation of historically marginalized demographic groups at all levels results in an unjust power imbalance (see the brief literature review in Section \ref{sec:survey} for further details).
Because of the nature of the dataset \cite{wapman-data}, which contains gender data but not other demographic data (and in particular no data on race or ethnicity), we are only able to conduct analyses using gender and not other identity characteristics.  
However, by combining this dataset with publicly available data from the National Science Foundation (NSF) on awards made by the Division of Mathematical Sciences (DMS) \cite{NSFDMS}, which is the primary funder of mathematical sciences research in the United States, we can investigate relationships among gender equity, \prestige, and NSF DMS funding.

Our thesis is that we can quantify and document inequity in mathematical sciences PhD-granting departments in the United States. 
Specifically, in this paper we will address the following research questions:

\begin{itemize}
    \item RQ1: What is the relationship between the percentage of women in a department and that department's \prestige?
    \item RQ2: What is the relationship between \prestige\ and funding received from the National Science Foundation Division of Mathematical Sciences?
    \item RQ3: What is the relationship between the percentage of women in a Mathematics department and funding received from the National Science Foundation Division of Mathematical Sciences?
\end{itemize}

\noindent The rest of this paper is organized as follows.
In Section \ref{sec:survey}, we provide some examples of recent research that uses data science and mathematics to analyze the mathematics and science communities. 
We describe the data and methods used in our research in Section \ref{sec:method}.  
Specifically, in this section we describe the processing of the Wapman et al. data and the NSF funding data that is required so that  we can investigate our research questions that support the thesis of this paper, i.e. that inequity exists in PhD-granting mathematical sciences in the United States.
We  provide  details about the statistical data analysis performed to establish the existence of quantifiable relationships between our variables of interest, gender percentage of faculty in mathematical sciences PhD-granting departments, amount of funding received from NSF from 2011-2020, and departmental \prestige.
The results of the data analysis and discussion of these results are given in Section \ref{sec:results} and Section \ref{sec:discussion}, respectively.  
We end the paper by discussing  in Section \ref{sec:future} some limitations of the work presented here and recognizing that there is a lot more work that can (and should)  be done to quantify and document inequity in the mathematical sciences.

\section{Existing work on inequity and hierarchy in mathematics}
\label{sec:survey}

In this section, we provide a short survey of selected recent work that uses tools, topics, and techniques from mathematics and data science to describe, document, and discuss inequity in the mathematical sciences.
We organize our summary of the literature in this area into three topics: 1) analysis of the (lack of) diversity in the mathematical sciences; 2)  existence of hierarchy in mathematics and other academic disciplines; and 3) evidence of inequity in the mathematical sciences. 
For a longer survey of the areas discussed here, as well as the broader field of the mathematics of Mathematics, we refer the reader to the recent paper by Buckmire et al. \cite{MathofMathSurvey}.

\subsection{The Demographics and Diversity   of the Mathematics Community.}
\label{DDMC}

It is well-documented that women are underrepresented at all levels in the mathematics community, and that their representation declines  as they progress through the academic system \cite{nasem2020women}. 
Women made up approximately 42.6\% of recipients of bachelor's degrees in mathematics between 2013 and 2018 \cite{buckmire2021diversity}. 
However,  only 29\% of doctorate recipients in mathematical sciences were women (in 2017-2018) \cite{AMSDoctorates}.  
Women were 28\% of hires in doctoral-granting mathematical sciences departments in 2019 \cite{AMSHiring}. 

Underrepresentation in the mathematical sciences based on race and ethnicity is profound and persistent \cite{ncses2023wmpd}. 
For example, from 2013 to 2018, the racial and ethnic makeup of undergraduate math graduates ranged from 64.9\% White, 7.5\% Latino/Hispanic, and 5.0\% Black to 52.5\% White, 9.9\% Latino/Hispanic and 4.2\% Black \cite{buckmire2021diversity}.
 Between 1993 and 2002, less than 5\% of those who earned doctoral degrees in mathematics were Black, Latino/Hispanic, or Indigenous, even though those communities made up a quarter of the general population of the United States at that time \cite{medina2004doctorate}. 
 
Vitulli \cite{vitulli2018gender} examined the representation of women being hired by mathematics departments, based upon data from annual surveys conducted by the American Mathematical Society (AMS). 
Prior to 2012\footnotemark, the AMS reported this data by dividing departments into three Groups based on the reputational rankings in the 1995 (or previously, 1982) National Research Council report on doctoral departments \cite{NRC82,NRC95}. 
Group I contained the highest rated 25.9\% of the departments. 
Group II was the next highest 30.3\% while Group III contained the remaining departments. 
Vitulli found that from 1991-2011, 20.5\% of the faculty hired by Group I departments were women, while 26.3\% of the faculty hired by the remaining departments were women. 
We acknowledge that this is just a small sample of the research in the literature analyzing the demographics and diversity of the mathematical sciences community.

\footnotetext{The AMS changed how they report this data in 2012, as the newest National Research Council report no longer provided a total ordering of departments, instead reporting multiple measures for each department.}

\subsection{The Existence of Hierarchy in Mathematics and Science.}
Recently, researchers have used available data on faculty positions at institutions of higher education in the United States to document the existence of hierarchies in faculty hiring networks in academia. 
These hierarchical structures \cite{kawakatsu2021emergence} in mathematics and science demonstrate that some institutions have greater influence on faculty hiring than others \cite{myers2011mathematical}; this is the essential characteristic of what we call \prestige. 
Institutions nearer the top of the hierarchy that are ``more prestigious" (i.e., they have greater \prestige) are more likely to have graduating doctoral students go on to obtain faculty positions at institutions that are ``less prestigious" and lower in the hierarchy.
Clauset et al. \cite{clauset2015systematic} demonstrated the existence of hierarchy in faculty hiring in a study involving departments in computer science, business, and history. 
Wapman et al. \cite{wapman2022quantifying} expanded this analysis to cover 295,089 faculty in 10,612 departments at 368 PhD-granting institutions and all academic disciplines for the years 2011-2020. 
FitzGerald et al. \cite{fitzgerald_huang_leisman_topaz_2022} built upon Wapman et al.'s research and partially replicated their results by using data from the Mathematics Genealogy Project (MGP) to restrict their analysis to mathematics faculty.
These results confirm that hierarchies exist in faculty hiring networks in the mathematical sciences. 


\subsection{Evidence of inequity in Mathematics.} 
There are multiple research articles that use quantitative tools and techniques to analyze and highlight examples of  inequity in the mathematics community. 
Topaz et al. \cite{topaz2016gender} analyzed the editorial boards of 435 mathematical science journals and found that women accounted for a mere 8.9\% of editorial positions.
Editorial positions play important gatekeeping roles and represent status in the mathematics community, so the underrepresentation of women in this area demonstrates inequity based on gender in the area of power over knowledge production.
Brisbin and Whitcher \cite{brisbin2015women} found that women are underrepresented as authors among papers in the mathematical sciences uploaded to the arXiv preprint repository, and that there are certain subfields (particularly concentrated within ``pure" or theoretical mathematics) with even larger discrepancies. 
Researchers have analyzed data describing different aspects of academic activity and demonstrated myriad ways that gender can negatively mediate opportunity for advancement, participation, and achievement in science and mathematics
\cite{mihaljevic2016effect,rissler2020gender,schmaling2023gender,lerman2022gendered}.
Schlenker \cite{schlenker2020prestige} notes that fields with applications to the social or physical sciences such as numerical analysis, mathematical modeling, or statistics (i.e. fields seen as being in applied mathematics) seem to be viewed by some as having low status in the wider mathematical community.
A large study investigating class backgrounds in academia by Morgan et al. \cite{morgan2022socioeconomic} found that faculty are much more likely than the general population to have a parent with a PhD, with the effect being even more pronounced at  institutions in the top quintile of U.S. News and World Report rankings.




\section{Data and Methods}
\label{sec:method}
In this section, we will describe the data and explain the methodology used to obtain our results.
Our data and code is publicly available at \cite{data}.
We utilized two datasets, one sourced from  Wapman et al.  \cite{wapman-data} and the second from awards made by the Division of Mathematical Sciences (DMS) at the United States National Science Foundation (NSF) between 2011 and 2020 \cite{NSFAwardsSearch}.

The Wapman et al. dataset required some nuance to interpret.  
The dataset consisted of a census of tenured or tenure-track faculty employed at \textit{all} PhD-granting institutions in the United States from the years 2011–2020. 
Faculty were only included in this sample if they were employed in the majority of the years under review.  
This dataset is centered on departments, rather than on faculty, and these departments are each assigned to a field such as ``Mathematics.''  
In particular, a department may be accounted for in multiple fields; a ``Department of Mathematics and Statistics'' would have its faculty included twice, in the fields of ``Mathematics'' and ``Statistics.''

Because our goal is to quantify and document inequity in the mathematics community as a whole, we choose to define the mathematical sciences  as broadly as possible (see  \cite{buckmire2023definitions}).
This choice results in a reduction of Wapman et al.'s original dataset  of  295,089 faculty in all academic disciplines offering PhDs in the United States to  9,814 faculty that are distributed among the fields of Mathematics, Statistics, or Operations Research (Table \ref{tab:counts}). 
We adopt the convention of capitalizing these three terms when referring to these fields present in the data throughout the rest of this paper.
\begin{table}[ht]
\centering
    \caption{Faculty present in the Wapman et al. dataset in the fields of Mathematics, Statistics, and Operations Research.} \label{tab:counts}
    \begin{tabular}{lccc}
    Field & Departments & Faculty Members& Percentage of women \\ \hline
    Mathematics & 223 & 7328 & 16.8\% \\
    Statistics & 122 & 2576 & 20.9\% \\
    Operations Research & 51 & 1034 & 19.3\% \\ 
    \end{tabular}
\end{table}
In our analysis, we incorporate the department prestige rankings from Wapman et al. \cite{wapman-data} which we refer to as \prestige.
Because of the way these are computed (a department has more \prestige\ if its PhD-graduates are hired by departments that also have \prestige), some departments do not have a \prestige\ ranking.
This could happen if the department does not have a PhD program in one of the three fields (recall the Wapman dataset began with {\em institutions} that grant PhDs), or if none of its PhD-graduates were hired by departments with \prestige\ rankings. 
In Mathematics, for example, there are 223 departments listed, but only 161 of these have a \prestige\ ranking.
Departments without \prestige\ rankings were not included in the data analyses involving \prestige\ below, but were included in the analysis of NSF funding later in the paper.
  
Since the Group I departments---in the groupings used by the AMS historically to divide departments by perceived quality (see Section \ref{DDMC} above)---constituted about 25\% of departments, in our analysis we consider the upper quartile (in terms of \prestige\ rankings) as ``elite" departments and compare this group to the remaining 75\% of departments, which we refer to as ``non-elite."

As stated earlier in the paper, the Wapman et al. dataset does not have information on race or ethnicity but does include gender.
However, gender  is only included as a binary variable in their dataset.  
In fact, gender is self-reported for a small percentage (6\%) of faculty in their initial dataset. 
They then attempted to infer the gender of the remaining faculty based on their names and the use of software that claims to assign gender based on names relatively accurately; ultimately  binary gender was ascribed to a total of 85\% of the faculty listed in their dataset. 
We include only these faculty who were ascribed a binary gender in our analyses,which eliminates roughly 15\% of the total due to the inability to accurately ascribe a gender to these data entries.

We separately obtained data from the NSF's publicly available data on awards made by the Division of Mathematical Sciences (DMS) in the Directorate of Mathematical and Physical Sciences (MPS) from 2011-2020 \cite{NSFAwardsSearch}. 
These data were aggregated by institution. 
Since the institution names in the Wapman et al. dataset typically did not match the formal organization names listed by the NSF, we manually adjusted these in order to compare the two datasets. 
On average, DMS awarded \$235 million per year towards achieving its mission to support ``a wide range of research in mathematics and statistics aimed at developing and exploring the properties and applications of mathematical structures'' \cite{NSFDMS}.
Of the awards to institutions, 80\% were matched to the departments of interest in the following analysis.
The NSF, and DMS in particular, is the primary source of funding for mathematics research in the United States \cite{NSFDMS}.  
While some mathematicians (such as several authors of this paper) receive funding from  NSF divisions other than DMS and directorates other than MPS (e.g., the Division of Undergraduate Education and the Directorate for STEM Education) as well as from other federal agencies (e.g., the National Institutes of Health), we consider DMS funding a reasonable measure for overall financial support of mathematics by the federal government.

There are two major caveats to our use of the NSF DMS funding data.
First, NSF awards are made to institutions rather than specifically to departments, which is the unit of analysis provided by the Wapman et al. data.
Second, since faculty in Statistics and Operations Research typically have more varied sources of funding,  we only considered the field of Mathematics in our analyses of funding discussed below.

Some of the institutions in the NSF data with the largest amount of funding were not included in our analysis because they were mathematics research institutes, not PhD-granting institutions.  
For example, the Mathematical Sciences Research Institute (MSRI, now known as the Simons Laufer Mathematical Sciences Institute or SLMath) and the Institute for Advanced Study (IAS) were awarded funds as separate entities from their associated universities, the University of California at Berkeley and Princeton University, respectively. 
We chose to treat MSRI and IAS as separate, non-PhD-granting institutions and therefore they are not included in our analyses.

\section{Results} \label{sec:results}
\begin{figure}[t]
    \centering
    \includegraphics[width=0.75\textwidth]{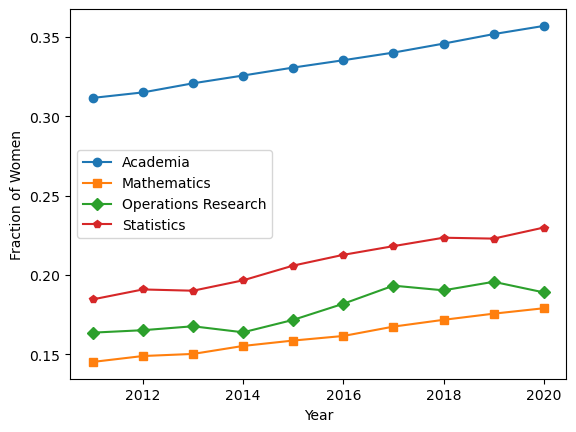}
    \caption{Fraction of women in the fields Mathematics, Statistics, and Operations Research, as well as academia as a whole, over the ten year period of the dataset. 
    }
    \label{fig:women-time}
\end{figure}

In this section, we present the primary results of our research into the distribution of faculty and funding at mathematical sciences PhD-granting institutions in the United States.
In Figure \ref{fig:women-time}, the percentage of women in the fields of Mathematics, Statistics, and Operations Research from 2011 to 2020 is given.
We note that the percentage of women in Mathematics lags behind Operations Research and Statistics throughout the time period of the dataset, mirroring the lower percentage of women that earn PhD's in Mathematics versus the other two fields \cite{AMSDoctorates}.
We further note the percentage of women in Mathematics, Statistics, and Operations Research is far below the percentage of women in Academia as a whole for the time period covered by the dataset. 

\begin{figure}[b]
\includegraphics[width=\textwidth]{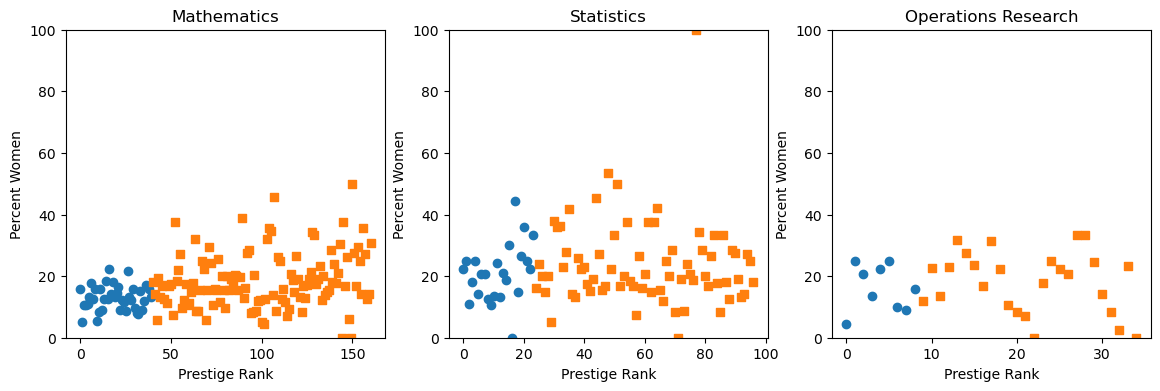}
\caption{Percentage of women by department in the fields of Mathematics, Statistics, and Operations Research. Color distinguishes the upper quartile of \prestige\ (blue circles) from the lower three quartiles (orange squares). 
} \label{fig:women-prestige}
\end{figure}

Next, we computed the percentages of faculty in each department inferred to be women, and plotted these according to \prestige\ rank  in the fields of Mathematics, Statistics, and Operations Research in Figure \ref{fig:women-prestige}. 
To compute this percentage, we used as a denominator the total number of faculty in a department for which a gender was inferred, in effect removing from our sample any faculty members whose gender could not be inferred. 
In Figure \ref{fig:women-prestige}, the blue circles are clustered in the lower-left corner of all three subfigures; this corresponds to the data demonstrating that elite departments (in the top quartile of \prestige) also have a low percentage of women (below 20\% for Mathematics).

In order to address RQ1, we then considered the elite  institutions as a group and calculated the percentage of faculty at these institutions that are women (Table \ref{tab:elite-women}); in each case, we see that the percentage of women among these elite institutions is lower than among non-elite institutions.  
A chi-squared test for each field was conducted, finding only the difference in Mathematics to be significant (\(p<0.001\)). 
We also conducted a Kendall tau test to determine if there is an association between the percentage of women and \prestige\ rank of Mathematics departments; we found a significant negative association between \prestige\ rank and rank by percentage of women (\(\tau=-0.23\), \(p<0.001\)). 
In other words, higher \prestige\ of Mathematics department is associated with having a lower percentage of women.


\begin{table}
\caption{The percentage of faculty at elite and non-elite institutions who are women in each field.}\label{tab:elite-women}
\begin{tabular}{lcc}
        Field & \multicolumn{1}{p{1.5in}}{Percentage of women among elite institutions} & \multicolumn{1}{p{1.7in}}{Percentage of women among non-elite institutions} \\ \hline
        Mathematics & 12.5\% & 18.1\% \\
        Statistics & 21.3\% & 21.6\% \\
        Operations Research &  17.0\% & 18.7\% \\
\end{tabular}
\end{table}

\begin{figure}[t!]
\centering
  \subcaptionbox{Total funding}[0.45\textwidth]{
\includegraphics[width=0.45\textwidth]{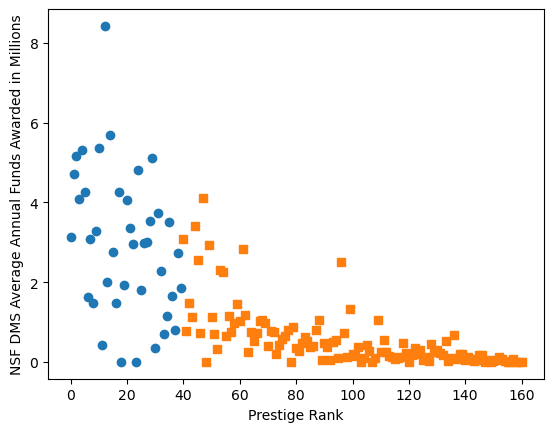} }
  \hfill
  \subcaptionbox{Funding per faculty}[0.45\textwidth]{
\includegraphics[width=0.45\textwidth]{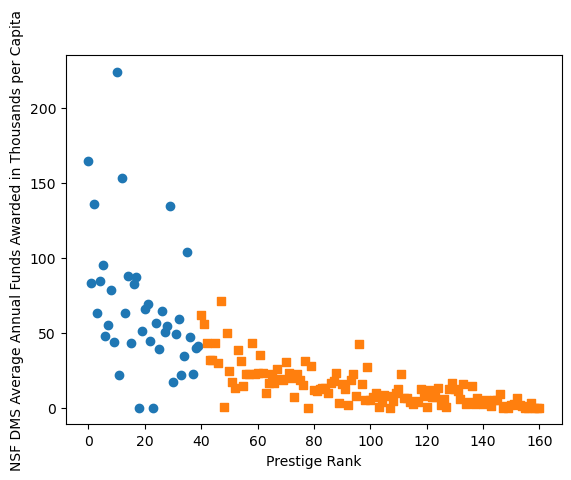}}
    \caption{Average annual grant funding for Mathematics departments by the prestige rank of the department, displayed by total funding to the department (a), and on a per capita basis accounting for the varying number of faculty in each department (b). 
    } \label{fig:funding-prestige}
\end{figure}

To address RQ2, we then explored the distribution of DMS funding to Mathematics departments with \prestige\ rankings from 2011-2020.
The elite institutions (i.e., the upper quartile by \prestige) were awarded in aggregate \$119M per year in grant funding, while the non-elite institutions, of which there are three times as many, were awarded only \$70M in aggregate of NSF money per year. 
We plotted this funding by \prestige\ of each department in Figure \ref{fig:funding-prestige}.
Note that the two subfigures in Figure~\ref{fig:funding-prestige} depict total funding and per capita funding versus \prestige\ and yet they are remarkably similar in appearance.
This demonstrates that we can use NSF funding to distinguish between different mathematical sciences PhD-granting institutions without worrying about individual outlier faculty members in Mathematics departments of those institutions.
We also computed the total amount of funding received by elite (top quartile by \prestige) and non-elite (lower 3 quartiles by \prestige) departments over the time period in question.
Elite departments received 64.7\% of funding in our dataset, compared to 35.3\% for the non-elite departments (despite these being thrice as numerous). 
Here we also conducted a Kendall tau test, finding a significant positive relationship between \prestige\ rank and DMS funding rank (\(\tau=0.68\), \(p<0.001\)). 
In other words, higher \prestige is associated with more DMS funding.

Another metric that can be used to quantify inequity in NSF DMS funding received by Mathematics departments in the United States  is the Gini coefficient, a well-known measure of inequality often used to characterize wealth inequality on a scale from 0 to 1. 
A uniform distribution of wealth would result in a Gini coeffient value of 0 and a case where one individual holds all the wealth would result in a Gini coefficient of 1. 
The Gini coefficient of NSF DMS funding for the Mathematics departments in our dataset is 0.63.
For more information on the Gini coefficient, refer to \cite{mukhopadhyay2021gini}. 

We also analyzed the full DMS-funded portfolio (excluding fellowships which are given to people instead of institutions) to avoid a possible sampling effect  due to our focus on  PhD-granting departments. 
In this more comprehensive dataset, the top 20\% holds 86.1\% of all DMS funding. 
The Gini coefficient of this distribution is 0.80. 
Thus, the larger set of all DMS funding recipients demonstrates greater inequality than the subset of PhD-granting institutions.

\begin{figure}[t]
\centering
    \includegraphics[width=0.75\textwidth]{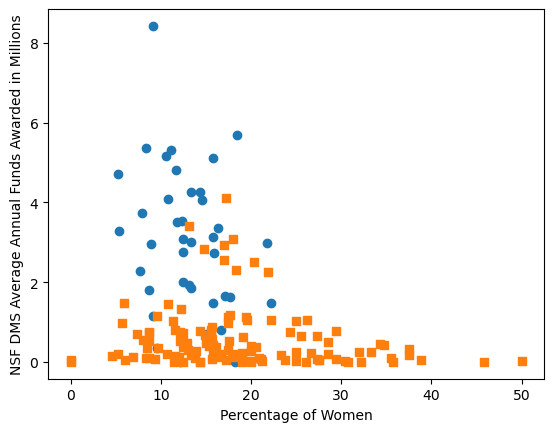}
    \caption{Average annual grant funding for Mathematics departments by the percentage of women faculty in the department. The upper quartile (by prestige) departments are represented by blue circles, while the remainder are orange squares. 
    } \label{fig:funding-women}
\end{figure}

To address RQ3, we  plotted the annual grant funding received by Mathematics departments against the percentage of women in those departments in Figure \ref{fig:funding-women}.  
We note an interesting effect in Figure~\ref{fig:funding-women}, where none of the 29 departments (which are all non-elite with respect to \prestige) with at least 25\% women received more than \$1.1M in average annual funding from the DMS. 
We conducted a Kendall tau test to determine if there is an association between annual funding received from DMS and the percentage of women in Mathematics departments and found a significant negative relationship (\(\tau=-0.22\), \(p<0.001\)).
In other words, a higher percentage of women is associated with less annual DMS funding,

In addition to the statistical analyses above, we also ran a k-means clustering analysis on the dataset that contained NSF DMS funding and faculty gender percentage for mathematical sciences PhD-granting  departments. 
We are particularly interested in the extent to which clusters produced by this  analysis would reflect the \prestige\ rankings.

\begin{figure}[h!]
    \centering
    \includegraphics[width=0.75\textwidth]{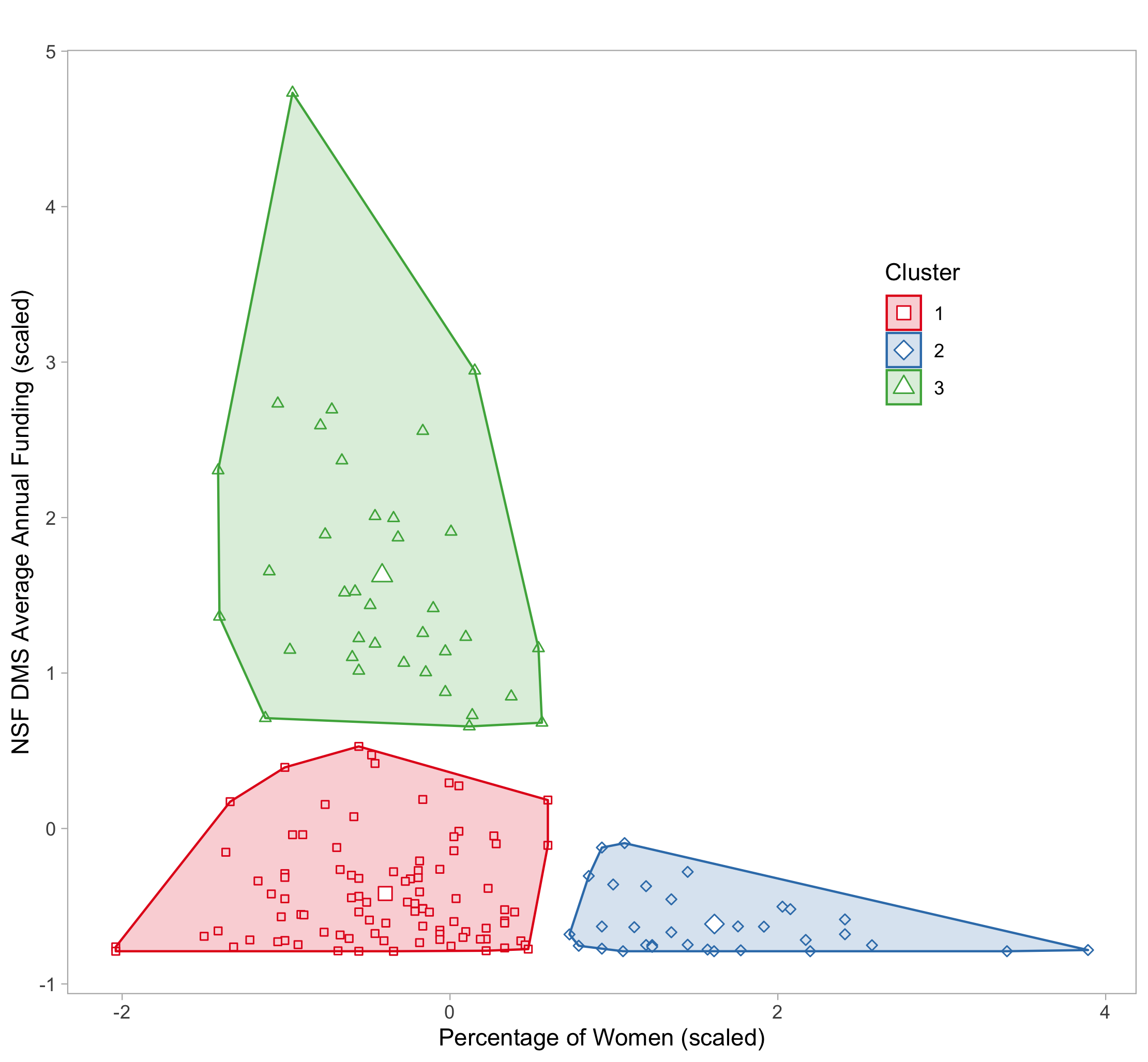}
    \caption{The clustering analysis of total NSF DMS funding amount versus the fraction of women by institution reveals three clusters. Cluster 1 in the bottom left, designated by red squares, identifies academic institutions with low representation of women and with low funding amounts. Cluster 2 in the bottom right, designated by blue diamonds,  identifies institutions  with high representation of women, but  low funding amounts. Cluster 3 in the top left, designated by green triangles, identifies institutions with low representation of women, but high funding amounts.} 
    \label{fig:cluster}
\end{figure}

After normalizing the data, we ran an analysis to determine the optimal clustering number. 
Using the \texttt{cluster} \cite{RCluster} and \texttt{factoextra} \cite{RFactoextra} packages in the R programming language, we identified three as a consistent optimal number in the ``bend-in-the-knee" analyses. 
These analyses identify the point at which increasing cluster number does not return significant advantages in lowering measures such as within-cluster variation. 
A visualization of the resulting three clusters is found in Figure \ref{fig:cluster}.
Cluster 1 in the bottom left, designated by red squares, identifies academic institutions with low representation of women and with low funding amounts.
Cluster 2 in the bottom right, designated by blue diamonds,  identifies institutions  with high representation of women, but with low funding amounts. 
Cluster 3 in the top left, designated by green triangles, identifies institutions with low representation of women, but high funding amounts. 

Note that there is no high gender diversity with high funding cluster, nor are there any institutions in the upper right area of Figure \ref{fig:cluster}.  In Table \ref{tab:cluster} the number of elite and non-elite institutions in each cluster are aggregated.
We then conducted a chi-squared test to see if cluster membership and ``eliteness" (being in the top 25\% of \prestige\ rankings) were independent; they were not (\(p<0.001\)).
In fact, Cluster 2 (blue diamonds) did not contain a single elite department, while 75\% of Cluster 3 (green triangles) were elite departments.   
The red cluster of institutions with low gender diversity and low funding includes a mix of elite and non-elite institutions. 
Every institution identified in the gender diverse blue cluster was also ascribed to the lower 75\% of the \prestige\ ranking list. 
The vast majority of institutions in the high funding green cluster were elite (in the top quartile by \prestige\ ranking). 
This cluster identified 68\% of upper quartile (elite) institutions and only 7\% of the lower three quartile (non-elite) institutions. 

\begin{centering}
\begin{table} 
\caption{The distribution of elite and non-elite Departments among the clusters identified using k-means clustering (see Figure~\ref{fig:cluster}).}\label{tab:cluster}
    \begin{tabular}{lcc}
    Cluster Name & Elite Departments & Non-elite Departments \\
    \hline 
    Cluster 1 (red squares) & 13 & 80 \\
    Cluster 2 (blue diamonds) & 0 & 32 \\
    Cluster 3 (green triangles) & 27  & 9  
    \end{tabular}
\end{table}
\end{centering}

\section{Discussion} \label{sec:discussion}
In this section we shall discuss the results presented above that demonstrate inequalities exist in faculty composition with respect to gender and the distribution of federal funding to mathematical sciences PhD-granting institutions in the United States.

The data shows that almost all PhD-granting institutions have Mathematics departments which are composed of faculty that are disproportionately male. 
In fact, not a single Mathematics department represented in this dataset was majority women. 
We found that the underrepresentation of women is more pronounced among elite Mathematics departments 
(recall that we defined ``elite'' departments as those in the upper quartile of departments in the prestige ranking generated by Wapman et al.).

We believe in the fundamental principle that mathematical talent is distributed equally among all groups of people who do mathematics. 
In the context of this paper, we therefore assume an equal distribution of mathematical talent among men and women. 
Under that assumption, the results presented here highlight the extreme level of underrepresentation of women in PhD-granting mathematical sciences departments.

Our analysis of NSF DMS funding identifies inequality in the amount of financial support for  mathematics PhD-granting departments depending on \prestige. 
Pareto models, also popularized as the ``20/80'' economic model, predict that approximately 80 percent of assets are held, gained or earned by only 20 percent of the population being studied \cite{pareto1964cours}.
We found that in the elite institutions, the top 25\% in our dataset by \prestige\ ranking, garnered 65\% of the total funds given to the subset of PhD-granting institutions with a \prestige\ ranking. 
When we examine all NSF DMS funding, the top 20\% of awardees receive 86\% of all funds, with a Gini coefficient of 0.8.
This result demonstrates a larger inequality than the classic ``20/80'' proportion. 

The procedures and policies that NSF uses to determine which institutions receive funding may actually reinforce  inequity in the mathematical sciences.
NSF uses a process called merit review, where anonymous reviewers are asked to read, rate, and review funding proposals submitted to the agency.
They are instructed to assess the intellectual merit and broader impacts of all proposals under five elements \cite{PAPPG}.  
Two of these elements in particular may bolster inequity in the mathematical sciences.
First, reviewers are asked ``How well qualified is the individual, team, or organization?''
This question is likely to skew reviewers towards considering institutional reputations since the information provided about qualifications typically includes individuals' institutional affiliation(s).  
Second, reviewers are asked ``Are there adequate resources available to the PI (either at the home organization or through collaborations) to carry out the proposed activities?''
This second question is likely to skew  reviewers to more positively rate proposals from  well-resourced institutions, which contributes to a rich-get-richer phenomenon that could explain the extremely inequitable, ``Pareto-like" distribution we found in our analysis of  NSF DMS funding data.

Even more troubling are the systemic effects that are propagated over time by this  inequality in the distribution of NSF DMS funding.
The clich\'{e} ``the rich-get-richer'' is a colloquial distillation of how systems that disproportionately allocate resources then have an easier time justifying disproportionate funding.
For example, well-resourced universities have significant funds internally for pilot projects.
They have teams of grant departments that assist in the writing and administration of grants, funded by high indirect cost rates. 
They also have research support teams devoted to data gathering and processing  as well as communications teams devoted to disseminating the results. 
We also note that expectations in obtaining grant funds vary widely between departments and institutions, and are often higher at the elite institutions (those in the top quartile using \prestige) in our dataset. 
This likely affects the number of proposals that are submitted by different institutions.
In short, the effect we are seeing in NSF DMS funding is likely a result of a complicated collection of processes which reinforce and exacerbate the status quo that is a demonstrably unequal system, despite often being presented as a meritocracy.

To summarize,  we have analyzed data on faculty in mathematical sciences PhD-granting departments and NSF DMS funding to institutions with these departments.
We found that higher \prestige\ is associated with a lower percentage of women in Mathematics departments (RQ1).   
We also found, perhaps unsurprisingly, that NSF DMS funding was associated with higher \prestige\ (RQ2).  
Finally, we found that NSF DMS funding to departments was associated with a lower percentage of women (RQ3).

In this paper, we have shown that being an elite institution is an indicator of a disproportionately high allocation of financial resources from the federal government and a low proportion of women; this documents the existence of systemic inequality, i.e inequity, in the mathematical sciences community.
We advocate for the redefinition of ``prestige'' (separate from \prestige) in mathematics in a way that better reflects equity for excellence \cite{NSFEES}. 
For example, one could define institutions that are more prestigious to be institutions that ``reflect the diversity of the US population'' \cite{NSFEES} among their faculty and doctoral graduates in mathematics.
Using this re-envisioned idea of prestige, we present the ten most prestigious PhD-granting Mathematics departments based on representation of women in Table \ref{tab:women-top-ten}. 
Here, we intentionally choose to turn away from \prestige, but provide an alternative definition of prestige based on gender representation.
One could also define ``prestige''  in the contribution that one makes not just to mathematics or academia, but in service to society.
This is a cultural shift, but one that could have significant advantages to producing a vision for science and mathematics that improves all lives.

\begin{table}
\centering
\caption{The top ten PhD-granting Mathematics departments by percentage of women.}\label{tab:women-top-ten}
    \begin{tabular}{lc}
    Institution & Percentage of Women \\ \hline
    Bryn Mawr College & 50.0\% \\
    Louisiana Tech University & 50.0\% \\
    University of California  Merced & 50.0\% \\
    Teachers College Columbia & 45.8\% \\
    University of Texas at Tyler &  40.0\% \\
    University of New Hampshire & 38.9\% \\
    Cleveland State University & 38.1\% \\
    Drew University & 37.5\% \\
    Illinois State University & 37.5\% \\
    Case Western Reserve University & 37.5\%
    \end{tabular}
    
\end{table}

\subsection{Limitations.}
There are a number of limitations that accompany the research presented in this paper that we want to highlight below.
The primary limitation is that the there is a paucity of publicly available, comprehensive, self-reported demographic data about the mathematics community.
The data we used was part of a dataset shared publicly by Wapman et al. \cite{wapman-data} after they had processed it, and the raw data was not available to us.
Their methodology of determining \prestige\  means that only PhD-granting departments are represented in the prestige data; a large number of faculty in the mathematical sciences who are at community colleges and predominantly undergraduate institutions are not included. 
Additionally, because of the way disciplines in the mathematical sciences are defined in the Wapman et al. dataset, faculty who are in ``Mathematics and Statistics" departments are counted in ``Mathematics" and again in ``Statistics."

Another important limitation is the way gender is ascribed to individuals in the Wapman et al. dataset. Recall that this dataset included gender, but only a small percentage were determined by the individuals themselves, with the rest inferred based on names.
While we acknowledge this practice is common, there are multiple issues with name-based gender inference.  
There is some degree of selection bias in which names can be ascribed to a particular (binary) gender; we note, again, that in the dataset used here, 15\% of entries were omitted due to an inability to confidently ascribe them a gender. 
Furthermore, the reduction of gender to a binary erases the experiences of gender-diverse mathematical scientists from this work.  

We also lament the lack of race/ethnicity in the data; there are important questions to be answered about the interaction of race/ethnicity and inequity in the mathematical sciences. 
A similar argument can be made about the absence of available data about LGBTQ+ identity.
The Wapman et al. dataset is limited to only tenured or tenure-track faculty.  
Without comprehensive demographic data, an intersectional analysis involving multiple identity characteristics is not possible.
As discussed below, we hope other researchers will collect or generate additional data that can be used to address important outstanding questions about the mathematical sciences discipline.

\section{Future Directions.}
\label{sec:future}
There are many other directions in which the research presented here could be extended.
It is important to study the questions addressed in this paper with respect to other dimensions of diversity, particularly marginalized social identities including race/ethnicity, sexual orientation, national origin, and disability status, among others.
This future work should be done in a way that allows analysis using intersections of multiple identity characteristics.

The addition of geographic location to the analysis of the gender diversity of PhD-granting institutions as well as of the distribution of federal funding is a possible direction of future research.

Another future direction is to expand this work to a wider range of institutions and faculty appointments.  
A study encompassing all types of institutions, and particularly community colleges, minority-serving institutions, and primarily undergraduate institutions, is necessary. 
Additionally, future work should investigate related questions about all types of faculty employed at these institutions, especially the increasing percentage of non-tenure track faculty.

This work's primarily goal has been to document inequity in the mathematical sciences. 
Future work could involve developing mathematical models that are informed by the data we have that documents inequity in the mathematical sciences and lead to a better understanding of the mechanisms involved.  
Specifically, the data \cite{buckmire2021diversity, AMSHiring, AMSDoctorates} showing the underrepresentation of women in new hires versus their underrepresentation in tenure-stream positions in PhD-granting institution raises some interesting questions that could also be addressed in future research.

We conclude by inviting interested researchers to join us in the ongoing MetaMath project to use mathematics and data science to analyze our discipline in order to promote social justice and enhance equity within all fields of the mathematical sciences. 

\subsection{Author Contributions.} All the listed co-authors contributed equally to the creation of this article; the order of attribution is alphabetical, as is customary in mathematics and is not intended to demonstrate any distinction in credit or effort.

\subsection*{Acknowledgements.}
This material is based upon work supported by the National Science Foundation under Grant No. DMS-1929284 while all authors were in residence at the Institute for Computational and Experimental Research in Mathematics in Providence, RI, during the Data Science and Social Justice: Networks, Policy, and Education program.
The content is solely the responsibility of the authors and does not necessarily represent the official views of the National Science Foundation.
The authors would like to thank Phil Chodrow and Victor Piercey for fruitful conversations that pushed this work forward, and Sam Zhang for contributions to a first draft of this work.
 
Buckmire acknowledges sabbatical support provided by the Office of the Dean of the College at Occidental College in Los Angeles, California.  

Diaz Eaton was supported in part by the Bates Enhanced Sabbatical Fund and Faculty Professional Development Fund.

Hibdon was supported, in part, by the National Institutes of Health’s National Cancer Institute, Grant Numbers U54CA202995, U54CA202997, and U54CA203000. The content is solely the responsibility of the authors and does not necessarily represent the official views of the National Institutes of Health.

Kauba is partially supported by the National Science Foundation South Carolina LSAMP
Bridge to Doctorate Fellowship HRD-2005030.

Pab\'on is partially supported by the National Science Foundation under Awards DMS-2108839 and DMS-1450182.

Vindas-Mel\'endez is partially supported by the National Science Foundation under Award DMS-2102921.



\bibliographystyle{plain}
\bibliography{bibliography}
\begin{description}
\item[Ron Buckmire] received his PhD in mathematics from Rensselaer Polytechnic Institute in 1994 and has been on the faculty of Occidental College in Los Angeles, California ever since. He is a passionate advocate for broadening the participation of people from historically marginalized and currently underrepresented groups in the mathematical sciences. 
Dr. Buckmire believes that: 1) mathematics is a human endeavor; 2) mathematics is created, discovered, learned, researched, and taught by people; 3) the identities and experiences of ``who does math" are important and that everyone and anyone can be and should be welcome in the mathematics community.

Mathematics Department, Occidental College, Los Angeles, CA 90032\\
ron@oxy.edu

\item[Carrie Diaz Eaton] has M.A. and PhD degrees in mathematics from University of Maine and University of Tennessee. She is now an associate professor of digital and computational studies, where her mixed methods teaching and research focus on social justice, complex adaptive systems and leadership, and their intersection with  society, STEM culture, and higher education. They are committed to fostering a new STEM ecosystem in which the current marginalized and oppressed are instead the liberated. 

RIOS Institute, Digital \& Computational Studies, Bates College, Lewiston, ME 04240\\ 
cdeaton@bates.edu

\item[Joseph E. Hibdon, Jr.] received his PhD in Applied Mathematics from Northwestern University in 2011.  He is an associate professor of mathematics at Northeastern Illinois University, a Hispanic Serving Institutions in Chicago, Illinois.  Dr. Hibdon works in many interdisciplinary teams across the country to increase diversity in STEM and broadening the participation of students in research at the undergraduate level.  Dr. Hibdon's recent research is in mathematical modeling, with a focus on public health and biological phenomena.

Department of Mathematics, Northeastern Illinois University, Chicago, IL 60625\\
j-hibdonjr@neiu.edu

\item[Jakini Kauba] is a PhD student in mathematical sciences at Clemson University. Her interests lie in stochastics and data science with an emphasis in interdisciplinary social justice work that addresses the many inequities of our nation. She believes in using mathematics to better serve underrepresented and marginalized communities. 

    School of Mathematical and Statistical Sciences, Clemson University, Clemson, SC, 29634 \\
    jkauba@g.clemson.edu

\item[Drew Lewis] received his PhD in mathematics from Washington University in St. Louis in 2012 and most recently served on the faculty at the University of South Alabama. He now works primarily in education research and faculty development.
drew.lewis@gmail.com

\item[Omayra Ortega] received her PhD in Applied Mathematics \& Computational Sciences from the University of Iowa in 2008. She is an Associate Professor of Applied Mathematics \& Statistics and the Assistant Dean of Research and Internships in the School of Science and Technology at Sonoma State University. Using the tools from statistics, mathematics, data science, public health and epidemiology, Dr. Ortega tackles the emerging health issues. She is deeply committed to broadening the participation of underrepresented minorities in STEM and mentoring students through the challenges of academia.

Department of Mathematics \& Statistics, Sonoma State University, Rohnert Park, CA 94928\\ 
ortegao@sonoma.edu

\item[José L. Pabón] is a PhD candidate in mathematics at the New Jersey Institute of Technology in Newark, N.J. He received his A.B. degree in mathematics from Princeton University in 2019 and has a particular affinity for working with disenfranchised populations, having first-hand experience with these given he was born and raised in Puerto Rico.

Department of Mathematical Sciences, Newark, New Jersey, 07102\\
jlp43@njit.edu

\item[Rachel Roca] is a PhD student in computational mathematics, science, and engineering at Michigan State University. Her interests lie in topological data analysis and computing education with an emphasis in social good. She believes in leveraging mathematical and data science tools for both advocacy and activism.

Department of Computational Mathematics, Science, and Engineering, Michigan State University, East Lansing, MI 48824\\
rocarach@msu.edu

\item[Andr\'es R. Vindas-Mel\'endez] received their PhD in mathematics from the University of Kentucky in 2021. 
Andr\'es is an NSF Postdoctoral Fellow and Lecturer at UC Berkeley and will begin as a tenure-track Assistant Professor at Harvey Mudd College in July 2024. 
Andr\'es' research interests are in algebraic, enumerative, and geometric combinatorics and have expanded to include applications of data science and mathematics for social justice. 

Department of Mathematics, University of California, Berkeley, CA 94720\\
Department of Mathematics, Harvey Mudd College, Claremont, CA 91711\\
andres.vindas@berkeley.edu; arvm@hmc.edu

\end{description}
\vfill\eject

\end{document}